\documentclass[twocolumn]{autart}    

\usepackage{graphicx}          
                               
\usepackage{amsmath}
\usepackage{amssymb}

\newcommand*\diff{\mathop{}\!\mathrm{d}}%

\usepackage{graphics} 
\usepackage{epsfig} 
\usepackage{epstopdf}
\usepackage{subfigure}
\usepackage{color}

\graphicspath{{./Figures/}}
\allowdisplaybreaks[4]

\begin{document}

\begin{frontmatter}

\title{An LMI Condition for the Robustness of Constant-Delay Linear Predictor Feedback with Respect to Uncertain Time-Varying Input Delays\thanksref{footnoteinfo}} 

\thanks[footnoteinfo]{This publication has emanated from research supported in part by a research grant from Science Foundation Ireland (SFI) under grant number 16/RC/3872 and is co-funded under the European Regional Development Fund and by I-Form industry partners.\\ 
Corresponding author H.~Lhachemi.}

\author[UCD]{Hugo Lhachemi}\ead{hugo.lhachemi@ucd.ie},
\author[GIPSA-lab]{Christophe Prieur}\ead{christophe.prieur@gipsa-lab.fr}, 
\author[UCD,Imperial]{Robert Shorten}\ead{robert.shorten@ucd.ie},               

\address[UCD]{School of Electrical and Electronic Engineering, University College Dublin, Dublin, Ireland}  
\address[GIPSA-lab]{Universit{\'e} Grenoble Alpes, CNRS, Grenoble-INP, GIPSA-lab, F-38000, Grenoble, France}             
\address[Imperial]{Imperial College London, South Kensington, UK}             

\begin{keyword}                           
Time-varying delay control; Predictor feedback; Robust stability; PDEs; Boundary control.             
\end{keyword}                             

\begin{abstract}                          
This paper discusses the robustness of the constant-delay predictor feedback in the case of an uncertain time-varying input delay. Specifically, we study the stability of the closed-loop system when the predictor feedback is designed based on the knowledge of the nominal value of the time-varying delay. By resorting to an adequate Lyapunov-Krasovskii functional, we derive an LMI-based sufficient condition ensuring the exponential stability of the closed-loop system for small enough variations of the time-varying delay around its nominal value. These results are extended to the feedback stabilization of a class of diagonal infinite-dimensional boundary control systems in the presence of a time-varying delay in the boundary control input.
\end{abstract}

\end{frontmatter}

\section{Introduction}
Originally motivated by the work of Artstein~\cite{artstein1982linear}, linear predictor feedback is an efficient tool for the feedback stabilization of Linear Time-Invariant (LTI) systems with constant input delay. In particular, predictor feedback can be used for controlling plants that are open-loop unstable and in the presence of large input delays. Many extensions have been reported (see, e.g., \cite{krstic2009delay} and the references therein). These include the case of time-varying delay linear predictor feedback~\cite{krstic2010lyapunov}; robustness with respect to disturbance signals~\cite{cai2018input}; truncated predictor~\cite{zhou2012truncated}; predictor observers in the case of sensor delays~\cite{krstic2009delay}; predictors for nonlinear systems~\cite{bekiaris2013robustness,krstic2010input}; dependence of the delay on the state~\cite{bekiaris2013nonlinear}; networked control~\cite{selivanov2016predictor}; and the boundary control of partial differential equations~\cite{lhachemi2018feedback,prieur2018feedback}.

Most of the predictor feedback strategies reported in the literature assume a perfect knowledge in real-time of the input delay. However, such an assumption might be difficult to fulfill in practice. Consequently, there has been an increased interest in the last decade for the study of the robustness of the predictor feedback with respect to delay mismatches. An example of such a problem was investigated in~\cite{krstic2008lyapunov} where the exponential stability of the closed-loop system was assessed for unknown constant delays with small enough deviations from the nominal value. The study of the impact of an unknown time-varying delay, but with known nominal value which is used to design the predictor feedback, on the system closed-loop stability was reported in~\cite{bekiaris2013robustness}. In particular, it was shown that the exponential stability of the closed-loop system is guaranteed for sufficiently small variations of the delay in both amplitude and rate of variation. Such an approach was further investigated in~\cite{karafyllis2013delay} where a small gain condition on the only amplitude of variation of the delay around its nominal value was derived for ensuring the exponential stability of the closed-loop system. However, as underlined in~\cite{li2014robustness}, such a small gain condition might be conservative as it involves norms of matrices which generally grow quickly with their dimensions. In order to reduce such a conservatism, it was proposed in~\cite{li2014robustness} to resort to a Lyapunov-Krasovskii functional approach in the case of constant uncertain delays. By doing so, an LMI-based sufficient condition, was derived, for ensuring the asymptotic stability of the closed-loop system with constant uncertain delays. 

The first contribution of this paper deals with the study of the robustness of the constant-delay predictor feedback that has been designed based on the nominal value of an uncertain and time-varying input delay. By taking advantage of classical Lyapunov-Krasovskii functionals~\cite{fridman2014tutorial}, we derive an LMI-based sufficient condition on the amplitude of variation of the input delay around its nominal value that ensures the exponential stability of the closed-loop system. Such an approch was investigated first in~\cite{selivanov2016predictor} in the context of networked control. However, the LMI condition derived in this paper differs from the one proposed in~\cite{selivanov2016predictor}. Three examples are developed showing that, for these case studies, the LMI condition proposed in this paper provides less conservative results than the small gain condition reported in~\cite{karafyllis2013delay} and the LMI condition extracted from~\cite{selivanov2016predictor}. 

The second contribution of this paper deals with the extension of the above result to the feedback stabilization of a class of diagonal infinite-dimensional boundary control systems~\cite{Curtain2012} in the presence of a time-varying delay in the boundary control input. The control strategy consists in 1) the use of a predictor feedback to stabilize a finite-dimensional subsystem capturing the unstable modes of the infinite-dimensional system; 2) ensuring that the control law designed on a finite-dimensional truncated models successfully stabilizes the full infinite-dimensional system. Such a control strategy, inspired by~\cite{russell1978controllability} in the case of a delay-free feedback control, was first reported in~\cite{prieur2018feedback} for the exponential stabilization of a reaction-diffusion equation with a constant delay in the boundary control. Note that a different approach for tackling the same feedback stabilization problem was reported in~\cite{krstic2009control} via the use of a backstepping boundary controller. Ideas from~\cite{prieur2018feedback} were extended to the exponential stabilization of a class of diagonal infinite-dimensional boundary control systems with constant delay in the boundary control in~\cite{lhachemi2018feedback}. In this present paper, we go beyond~\cite{krstic2009control,lhachemi2018feedback,prieur2018feedback} and assess the robustness of the control strategy reported in~\cite{lhachemi2018feedback} in the case of an uncertain and time-varying input delay. Specifically, we show that for time-varying delays presenting 1) a sufficiently small amplitude of variation around its nominal value (with sufficient condition provided by the LMI condition discussed above); 2) a rate of variation that is bounded by an arbitrarily large constant; the infinite-dimensional closed-loop system is exponentially stable.

The remainder of this paper is organized as follows. The robustness of the predictor feedback with respect to uncertain and time-varying delays is investigated in Section~\ref{sec: finite dimensional systems}. The extension of this result to the feedback stabilization of a class of diagonal infinite-dimensional boundary control systems is presented in Section~\ref{sec: infinite dimensional systems}. The obtained results are applied in Section~\ref{sec: illustrative example}. Finally, concluding remarks are provided in Section~\ref{sec: conclusion}.

\textbf{Notation.} The sets of non-negative integers, positive integers, real, non-negative real, positive real, and complex numbers are denoted by $\mathbb{N}$, $\mathbb{N}^*$, $\mathbb{R}$, $\mathbb{R}_+$, $\mathbb{R}_+^*$, and $\mathbb{C}$, respectively. The real and imaginary parts of a complex number $z$ are denoted by $\operatorname{Re} z$ and $\operatorname{Im} z$, respectively. The field $\mathbb{K}$ denotes either $\mathbb{R}$ or $\mathbb{C}$. The set of $n$-dimensional vectors over $\mathbb{K}$ is denoted by $\mathbb{K}^n$ and is endowed with the Euclidean norm $\Vert x \Vert = \sqrt{x^* x}$. The set of $n \times m$ matrices over $\mathbb{K}$ is denoted by $\mathbb{K}^{n \times m}$ and is endowed with the induced norm denoted by $\Vert\cdot\Vert$. For any symmetric matrix $P \in \mathbb{R}^{n \times n}$, $P \succ 0$ (resp. $P \succeq 0$) means that $P$ is positive definite (resp. positive semi-definite). The set of symmetric positive definite matrices of order $n$ is denoted by $\mathbb{S}_n^{+*}$. For any symmetric matrix $P \in \mathbb{R}^{n \times n}$, $\lambda_m(P)$ and $\lambda_M(P)$ denote the smallest and largest eigenvalues of $P$, respectively. For $M = (m_{i,j}) \in \mathbb{C}^{n \times m}$, we introduce
\begin{equation*}
\mathcal{R} (M) \triangleq 
\begin{bmatrix}
\operatorname{Re} M & -\operatorname{Im}M \\ \operatorname{Im} M & \operatorname{Re} M
\end{bmatrix} 
\in \mathbb{R}^{2n \times 2m}
\end{equation*}
where $\operatorname{Re} M \triangleq (\operatorname{Re} m_{i,j}) \in \mathbb{R}^{n \times m}$ and $\operatorname{Im} M \triangleq (\operatorname{Im} m_{i,j}) \in \mathbb{R}^{n \times m}$. For any $t_0>0$, we say that $\varphi \in \mathcal{C}^0(\mathbb{R};\mathbb{R})$ is a \textit{transition signal over $[0,t_0]$} if $0 \leq \varphi \leq 1$, $\left. \varphi \right\vert_{(-\infty,0]} = 0$, and $\left. \varphi \right\vert_{[t_0,+\infty)} = 1$. In Section~\ref{sec: infinite dimensional systems}, the notations and terminologies for infinite-dimensional systems are retrieved from~\cite{Curtain2012}.

\section{Delay-robustness of predictor feedback for LTI systems}\label{sec: finite dimensional systems}

\subsection{Problem setting and existing result}

The first part of this paper deals with the feedback stabilization of the following LTI system with delay control:
\begin{equation}\label{eq: LTI system with time-varying delay}
\dot{x}(t) = A x(t) + B u(t-D(t)) , \quad t \geq 0 ,
\end{equation}
with $A \in \mathbb{R}^{n \times n}$ and $B \in \mathbb{R}^{n \times m}$ such that the pair $(A,B)$ is stabilizable. Vectors $x(t) \in \mathbb{R}^n$ and $u(t) \in \mathbb{R}^m$ denote the state and the control input, respectively. The command input is subject to an uncertain time-varying delay $D \in \mathcal{C}^0(\mathbb{R}_+;\mathbb{R}_+)$. We assume that there exist $D_0 > 0$ and $0 < \delta < D_0$ such that $\vert D(t) - D_0 \vert \leq \delta$ for all $t \geq 0$. In this context, the following constant-delay linear predictive feedback, which is based on the knowledge of the constant nominal value $D_0$, has been proposed in~\cite{bekiaris2013robustness}:
\begin{equation}\label{eq: nominal predictor feedback}
u(t) = K \left\{ e^{D_0 A} x(t) + \int_{t-D_0}^{t} e^{(t-s)A} B u(s) \diff s \right\} 
\end{equation}
for $t \geq 0$, where $K \in \mathbb{R}^{n \times m}$ is a feedback gain such that $A_\mathrm{cl} \triangleq A+BK$ is Hurwitz. The validity of such a control strategy was assessed in~\cite{karafyllis2013delay} via a small gain argument.

\begin{thm}[\cite{karafyllis2013delay}]\label{thm: exp stab - small gain}
Let $D_0 > 0$ be given and let $A \in \mathbb{R}^{n \times n}$, $B \in \mathbb{R}^{n \times m}$, and $K \in \mathbb{R}^{m \times n}$ be such that $A_\mathrm{cl} = A+BK$ is Hurwitz. Let $\delta \in (0,D_0)$ be such that
\begin{equation}\label{eq: thm exp stab - small gain condition}
M \Vert e^{D_0 A}BK \Vert \left\{ e^{\Vert A_\mathrm{cl} \Vert \delta} - e^{-\mu \delta} \right\} < \mu ,
\end{equation}  
where $M,\mu > 0$ are constants satisfying $\Vert e^{A_\mathrm{cl} t} \Vert \leq M e^{- \mu t}$ for all $t \geq 0$. Then, there exists $N,\sigma > 0$ such that for all $x_0 \in \mathbb{R}^n$, $u_0 \in \mathcal{C}^0([-D_0-\delta,0];\mathbb{R}^m)$ with $u_0(0) = K \left\{ e^{D_0 A} x_0 + \int_{-D_0}^{0} e^{-sA} B u_0(s) \diff s \right\}$, and $D \in \mathcal{C}^0(\mathbb{R}_+;\mathbb{R}_+)$ with $\vert D - D_0 \vert \leq \delta$, the solution of (\ref{eq: LTI system with time-varying delay}-\ref{eq: nominal predictor feedback}) associated with the initial conditions $x(0) = x_0$ and $u(t) = u_0(t)$ for $-D_0-\delta \leq t \leq 0$ satisfies for all $t \geq 0$ the following estimate:
\begin{align*}
& \Vert x(t) \Vert + \max\limits_{t-D_0-\delta \leq s \leq t} \Vert u(s) \Vert  \\
& \qquad \leq N e^{-\sigma t} \left\{ \Vert x_0 \Vert + \max\limits_{-D_0-\delta \leq s \leq 0} \Vert u_0(s) \Vert \right\} .
\end{align*}
\end{thm}

As the left hand-side of (\ref{eq: thm exp stab - small gain condition}) is equal to zero when $\delta = 0$, a continuity argument shows that there always exists a $\delta > 0$ such that (\ref{eq: thm exp stab - small gain condition}) holds true. Therefore, Theorem~\ref{thm: exp stab - small gain} ensures the existence of a sufficiently small amplitude of perturbation $\delta > 0$ of the delay $D(t)$ around its nominal value $D_0$ such that the constant-delay linear predictor feedback (\ref{eq: nominal predictor feedback}) ensures the exponential stability of the closed-loop system with uncertain time-varying input delays. However, due to the nature of the small gain-condition (\ref{eq: thm exp stab - small gain condition}) that involves the norm of matrices (which generally grow quickly as a function of the matrices dimensions $n$ and $m$), the admissible values of $\delta$ might be conservative (see~\cite{li2014robustness}). In particular, from the fact that $M \geq 1$ and $0 < \mu \leq \mu_M(A_\mathrm{cl}) \triangleq - \max\{\operatorname{Re} \lambda \,:\, \lambda\in\mathrm{sp}_\mathbb{C}(A_\mathrm{cl})\}$, any $\delta > 0$ such that the small gain condition (\ref{eq: thm exp stab - small gain condition}) holds true satisfies the following estimate:
\begin{equation}\label{eq: thm exp stab - small gain condition - estimate}
\delta < \delta_\mathrm{E} \triangleq \dfrac{1}{\Vert A_\mathrm{cl} \Vert} \log\left( 1 + \dfrac{\mu_M(A_\mathrm{cl})}{\Vert e^{D_0 A} BK \Vert} \right) .
\end{equation}
To reduce the conservatism, an LMI condition ensuring the exponential stability of the closed-loop system was derived in~\cite{selivanov2016predictor} in the context of networked control. The objective of this section it to propose the construction of an alternative LMI for such a problem. Numerical comparisons between the different methods (small gain and LMIs) will be carried out in Subsection~\ref{subsec: applications} and Section~\ref{sec: illustrative example}.

\subsection{Preliminary results}

For $h > 0$, we denote by $W$ the space of absolutely continuous functions $\psi : [-h,0] \rightarrow \mathbb{R}^n$ with square-integrable derivative endowed with the norm $\Vert \psi \Vert_{W} \triangleq \sqrt{ \Vert \psi(0) \Vert^2 + \int_{-h}^0 \Vert \dot{\psi}(\theta) \Vert^2 \diff \theta}$ (see~\cite[Chap.~4, Sec.~1.3]{kolmanovskii2012applied}).

\begin{lem}\label{lem: lemma 1}
Let $M,N \in \mathbb{R}^{n \times n}$, $D_0 > 0$, and $\delta \in (0,D_0)$ be given. Assume that there exist $\kappa > 0$, $P_1,Q\in\mathbb{S}_n^{+*}$, and $P_2, P_3 \in \mathbb{R}^{n \times n}$ such that $\Theta(\delta,\kappa) \preceq 0$ with
\begin{align}
& \Theta(\delta,\kappa) = \label{eq: preleminary lemma LMI condition} \\
& \begin{bmatrix}
2 \kappa P_1 + M^\top P_2 + P_2^\top M & \;P_1 - P_2^\top + M^\top P_3 & \;\delta P_2^\top N \\
P_1 - P_2 + P_3^\top M & \; - P_3 - P_3^\top + 2 \delta Q & \;\delta P_3^\top N \\
\delta N^\top P_2 & \; \delta N^\top P_3 & \; - \delta e^{-2 \kappa D_0} Q
\end{bmatrix} . \nonumber 
\end{align}
Then, there exists $C_0>0$ such that, for any  $D \in \mathcal{C}^0(\mathbb{R}_+;\mathbb{R}_+)$ with $\vert D - D_0 \vert \leq \delta$, the trajectory $x$ of:
\begin{align*}
\dot{x}(t) & = M x(t) + N \left\{ x(t-D(t)) - x(t-D_0) \right\} , \quad t \geq 0 ; \\
x(t) & = x_0(t) , \quad t \in [-D_0-\delta,0] 
\end{align*}
with initial condition $x_0 \in W$ (for $h = D_0+\delta$) satisfies $\Vert x(t) \Vert \leq C_0 e^{- \kappa t} \Vert x_0 \Vert_W$ for all $t \geq 0$. 
\end{lem}

\textbf{Proof.} 
For all $t \geq 0$, one has
\begin{equation}
\dot{x}(t) = M x(t) + N \int_{t-D_0}^{t-D(t)} \dot{x}(\tau) \diff \tau . \label{eq: dynamics in function of dot_x}
\end{equation}
Inspired by classical Lyapunov-Krasovskii functional depending on time derivative for systems with fast varying delays, see~ \cite[Sec.~3.2]{fridman2014tutorial}, we introduce $V(t) = V_1(t) + V_2(t)$ with $V_1(t) = x(t)^\top P_1 x(t)$ and $V_2(t) = \int_{-D_0-\delta}^{-D_0+\delta} \int_{t+\theta}^{t} e^{2\kappa(s-t)} \dot{x}(s)^\top Q \dot{x}(s) \diff s \diff \theta$, where $P_1,Q\in\mathbb{S}_n^{+*}$. Then we have, for all $t \geq 0$,
\begin{align}
\dot{V}(t)
& = 2 x(t)^\top P_1 \dot{x}(t) + 2 \delta \dot{x}(t)^\top Q \dot{x}(t) - 2 \kappa V_2(t) \label{eq: computation dot_V} \\
& \phantom{=}\, - \int_{-D_0-\delta}^{-D_0+\delta} e^{2 \kappa \theta} \dot{x}(t+\theta)^\top Q \dot{x}(t+\theta) \diff \theta . \nonumber
\end{align}
The remaining of the proof is now an adaptation of~\cite[Proof of Thm~1]{fridman2006new}. Introducing $P = \begin{bmatrix} P_1 & 0 \\ P_2 & P_3 \end{bmatrix}$, where $P_2, P_3 \in \mathbb{R}^{n \times n}$ are ``slack variables''~\cite{fridman2014tutorial}, we have
\begin{align}
& x(t)^\top P_1 \dot{x}(t) \nonumber \\
& \overset{(\ref{eq: dynamics in function of dot_x})}{=}  
\begin{bmatrix}
x(t) \\ \dot{x}(t)
\end{bmatrix}^\top
P^\top
\begin{bmatrix}
\dot{x}(t) \\ -\dot{x}(t) + Mx(t) + N \int_{t-D_0}^{t-D(t)} \dot{x}(\tau) \diff \tau
\end{bmatrix} \nonumber  \\
& = 
\begin{bmatrix}
x(t) \\ \dot{x}(t)
\end{bmatrix}^\top
P^\top
\begin{bmatrix}
0 & I \\ M & -I 
\end{bmatrix}
\begin{bmatrix}
x(t) \\ \dot{x}(t)
\end{bmatrix} \label{eq: coumputation x'*P1*dot_x} \\
& \phantom{=}\; + \int_{t-D_0}^{t-D(t)}  
\begin{bmatrix}
x(t) \\ \dot{x}(t)
\end{bmatrix}^\top
P^\top
\begin{bmatrix}
0 \\ N
\end{bmatrix}
\dot{x}(\tau) \diff \tau . \nonumber 
\end{align}
Now, from the fact that, for any $a,b \in \mathbb{R}^n$, $2 a^\top b \leq \Vert a \Vert^2 + \Vert b \Vert^2$, we obtain that
\begin{align*}
& 2 \begin{bmatrix}
x(t) \\ \dot{x}(t)
\end{bmatrix}^\top
P^\top
\begin{bmatrix}
0 \\ N
\end{bmatrix}
\dot{x}(\tau) \\
& = 2 \left( e^{- \kappa (\tau-t)} Q^{-1/2} 
\begin{bmatrix}
0 \\ N
\end{bmatrix}^\top
P
\begin{bmatrix}
x(t) \\ \dot{x}(t)
\end{bmatrix} \right)^\top
\left( e^{\kappa (\tau-t)} Q^{1/2} \dot{x}(\tau) \right) \\
& \leq 
e^{-2 \kappa (\tau-t)}
\begin{bmatrix}
x(t) \\ \dot{x}(t)
\end{bmatrix}^\top
P^\top
\begin{bmatrix}
0 \\ N
\end{bmatrix}
Q^{-1}
\begin{bmatrix}
0 \\ N
\end{bmatrix}^\top
P
\begin{bmatrix}
x(t) \\ \dot{x}(t)
\end{bmatrix} \\ 
& \phantom{\leq}\,  + e^{2 \kappa (\tau-t)} \dot{x}(\tau)^\top Q \dot{x}(\tau) .
\end{align*}
With (\ref{eq: computation dot_V}-\ref{eq: coumputation x'*P1*dot_x}) we deduce that
\begin{align*}
& \dot{V}(t) + 2 \kappa V(t)\\
& \leq 2 \kappa V_1(t) +
2
\begin{bmatrix}
x(t) \\ \dot{x}(t)
\end{bmatrix}^\top
P^\top
\begin{bmatrix}
0 & I \\ M & -I 
\end{bmatrix}
\begin{bmatrix}
x(t) \\ \dot{x}(t)
\end{bmatrix} \\
& \phantom{\leq}\, + \delta e^{2 \kappa D_0}
\begin{bmatrix}
x(t) \\ \dot{x}(t)
\end{bmatrix}^\top
P^\top
\begin{bmatrix}
0 \\ N
\end{bmatrix}
Q^{-1}
\begin{bmatrix}
0 \\ N
\end{bmatrix}^\top
P
\begin{bmatrix}
x(t) \\ \dot{x}(t)
\end{bmatrix} \\
& \phantom{\leq}\, + 2 \delta \dot{x}(t)^\top Q \dot{x}(t)
+ \int_{t-D_0}^{t-D(t)} e^{2 \kappa (\tau-t)} \dot{x}(\tau)^\top Q \dot{x}(\tau) \diff \tau \\
& \phantom{\leq}\, - \int_{t-D_0-\delta}^{t-D_0+\delta} e^{2 \kappa (\tau-t)} \dot{x}(\tau)^\top Q \dot{x}(\tau) \diff \tau \\
& \leq 
\begin{bmatrix}
x(t) \\ \dot{x}(t)
\end{bmatrix}^\top
\left\{
\Psi
+ \delta e^{2 \kappa D_0}
P^\top
\begin{bmatrix}
0 \\ N
\end{bmatrix}
Q^{-1}
\begin{bmatrix}
0 \\ N
\end{bmatrix}^\top
P
\right\}
\begin{bmatrix}
x(t) \\ \dot{x}(t)
\end{bmatrix} ,
\end{align*}
where it has been used the fact that the sum of the two integral terms is always non positive, and with
\begin{equation*}
\Psi \triangleq
P^\top
\begin{bmatrix}
0 & I \\ M & -I
\end{bmatrix}
+
\begin{bmatrix}
0 & I \\ M & -I
\end{bmatrix}^\top
P
+ 2
\begin{bmatrix}
\kappa P_1 & 0 \\ 0 & \delta Q 
\end{bmatrix} .
\end{equation*}
From $\Theta(\delta,\kappa) \preceq 0$, the direct application of the Schur complement yields $\dot{V}(t) + 2 \kappa V(t) \leq 0$. The conclusion follows from the fact that $\lambda_\mathrm{m}(P_1) \Vert x(t) \Vert^2 \leq  V(t) \leq \max\left( \lambda_\mathrm{M}(P_1) , 2 \delta \lambda_\mathrm{M}(Q) \right) \Vert x(t+\cdot) \Vert_W^2$ for all $t \geq 0$ . \qed

By a continuity argument, $\Theta(\delta,0) \prec 0$ implies $\Theta(\delta,\kappa) \preceq 0$ for some $\kappa>0$. We deduce the following result.

\begin{cor}
Let $M,N \in \mathbb{R}^{n \times n}$, $D_0 > 0$, and $\delta \in (0, D_0)$ be given. Assume that $\Theta(\delta,0) \prec 0$. Then the conclusions of Lemma~\ref{lem: lemma 1} hold true for some decay rate $\kappa > 0$.
\end{cor}

From Lemma~\ref{lem: lemma 1}, the feasibility of the LMI $\Theta(\delta,\kappa) \preceq 0$ ensures that $M$ is Hurwitz. The following lemma states a form of converse result.

\begin{lem}\label{lem: lemma 2}
Let $D_0 > 0$ and $M,N \in \mathbb{R}^{n \times n}$ with $M$ Hurwitz be given. Let $P_2\in\mathbb{S}_n^{+*}$ be the unique solution of $M^\top P_2 + P_2 M = - I$ and let $0 \leq \kappa < 1/(4 \lambda_M(P_2))$ be given. Introducing $\delta^* = \delta^*(\kappa) > 0$ defined by\footnote{With the convention $\delta^* = D_0$ in the case $N = 0$.}
\begin{equation*}
\delta^* \triangleq
\min\left( 
D_0
,
\dfrac{\min\left\{ 1 - 4 \kappa \lambda_M(P_2) , \lambda_m\left( (M^{-1})^\top M^{-1} \right) \right\}}{2 \sqrt{2} e^{\kappa D_0} \left\Vert N^\top \begin{bmatrix} P_2 & -(M^{-1})^\top P_2  \end{bmatrix}  \right\Vert}
\right) ,
\end{equation*}
the LMI $\Theta(\delta,\kappa) \prec 0$ is feasible for all $\delta \in (0,\delta^*)$.
\end{lem}

\textbf{Proof.} 
As $M$ is Hurwitz, let $P_2\in\mathbb{S}_n^{+*}$ be the unique solution of the Lyapunov equation $M^\top P_2 + P_2 M = - I$. We introduce $P_1=2P_2\in\mathbb{S}_n^{+*}$, $P_3 = -(M^{-1})^\top P_2$, and $Q = \alpha I\in\mathbb{S}_n^{+*}$ with $\alpha > 0$. Then $\Theta(\delta,\kappa) \prec 0$ becomes:
\begin{equation}\label{eq: coro existence delta>0 - LMI condition}
\begin{bmatrix}
4 \kappa P_2 -I & 0 & \; \delta P_2 N \\
0 & - S_3 + 2 \alpha\delta I & \; \delta P_3^\top N \\
\delta N^\top P_2 & \delta N^\top P_3 & \; - \alpha\delta e^{-2 \kappa D_0} I
\end{bmatrix} 
\prec 0 .
\end{equation}
with $S_3 = P_3+P_3^\top = (M^{-1})^\top M^{-1} \succ 0$. As $\alpha,\delta>0$, the Schur complement shows that (\ref{eq: coro existence delta>0 - LMI condition}) is equivalent to
\begin{align*}
& \begin{bmatrix}
4 \kappa P_2 -I & 0 \\
0 & - S_3 + 2 \alpha\delta I \\
\end{bmatrix}
+ \dfrac{\delta}{\alpha} e^{2 \kappa D_0}
\begin{bmatrix}
P_2 N \\
P_3^\top N 
\end{bmatrix} 
\begin{bmatrix}
P_2 N \\
P_3^\top N 
\end{bmatrix}^\top
\prec 0 .
\end{align*}
A sufficient condition ensuring that the above LMI is satisfied is provided by $- \beta_0 + 2 \alpha \delta + \dfrac{\delta}{\alpha} \beta_1 < 0$, where $\beta_0 \triangleq \min\left\{ 1 - 4 \kappa \lambda_M(P_2) , \lambda_m (S_3) \right\} > 0$ and $\beta_1 \triangleq e^{2 \kappa D_0} \left\Vert N^\top \begin{bmatrix} P_2 & P_3  \end{bmatrix}  \right\Vert^2 \geq 0$. We deduce that $\delta <
\dfrac{\beta_0}{2\alpha + \beta_1/\alpha}$ implies that $\Theta(\delta,\kappa) \prec 0$, where $\alpha > 0$ can be freely selected. In the case $N = 0$, we obtain that $\beta_1 = 0$ and thus, by letting $\alpha \rightarrow 0^+$, $\delta^* = D_0$. In the case $N \neq 0$, we have $\beta_1 > 0$. Indeed, by contradiction, $\beta_1 = 0$ implies $N^\top P_2 = P_2 N = 0$. Multiplying $M^\top P_2 + P_2 M = - I$ from the left side by $N^\top$ and from the right side by $N$, we obtain that $N^\top N = 0$ yielding $N=0$. To conclude the proof, it is sufficient to note that, for any given $a,b>0$, the function $f(\alpha) = a \alpha + b/\alpha$ is such that $f(\alpha) \geq f(\sqrt{b/a})=2\sqrt{ab}$ for all $\alpha > 0$. \qed

\subsection{Robustness of constant-delay predictor feedback with respect to time-varying input delays}

We can now introduce the main result of this section.

\begin{thm}\label{thm: exp stability for finite-dimensional systems}
Let $A \in \mathbb{R}^{n \times n}$, $B \in \mathbb{R}^{n \times m}$, and $K \in \mathbb{R}^{m \times n}$ be such that $A_\mathrm{cl} \triangleq A + BK$ is Hurwitz. Let $\varphi$ be a transition signal\footnote{See notation section.} over $[0,t_0]$ with $t_0 > 0$ and let $D_0 > 0$ be a given nominal delay. Then, there exists $\delta \in (0,D_0)$ such that for any  $D \in \mathcal{C}^0(\mathbb{R}_+;\mathbb{R}_+)$ with $\vert D - D_0 \vert \leq \delta$, the closed-loop system given for $t \geq 0$ by 
\begin{align*}
\dot{x}(t) & = A x(t) + B u(t-D(t)) ,  \\
\left. u \right\vert_{[-D_0-\delta,0]} & = 0 , \\
u(t) & = \varphi(t) K e^{D_0 A} x(t) \\
& \phantom{=}\, + \varphi(t) K \int_{t-D_0}^t e^{(t-s)A}B u(s) \diff s , \nonumber \\
x(0) & = x_0  
\end{align*}
with initial condition $x_0 \in \mathbb{R}^n$ is exponentially stable in the sense that there exist constants $\kappa,C_1>0$, independent of $x_0$ and $D$, such that $\Vert x(t) \Vert + \Vert u(t) \Vert \leq C_1 e^{-\kappa t} \Vert x_0 \Vert$. In particular, this conclusion holds true (resp., with given decay rate $\kappa > 0$) for any $\delta \in (0,D_0)$ such that there exist $P_1,Q\in\mathbb{S}_n^{+*}$ and $P_2, P_3 \in \mathbb{R}^{n \times n}$ for which the LMI $\Theta(\delta,0) \prec 0$ (resp., $\Theta(\delta,\kappa) \preceq 0$) holds true with $M = A_\mathrm{cl}$ and $N = e^{D_0 A} BK$.
\end{thm}

\textbf{Proof.} 
Let $\delta \in (0,D_0)$ be such that $\Theta(\delta,0) \prec 0$ is feasible (see Lemma~\ref{lem: lemma 2}) and, by a continuity argument, let $\kappa > 0$ be such that $\Theta(\delta,\kappa) \preceq 0$ is feasible. By the properties of the Artstein transformation~\cite{bresch2018new}, we have $x \in \mathcal{C}^1(\mathbb{R}_+;\mathbb{R}^n)$ and $u \in \mathcal{C}^0([-D_0-\delta,+\infty);\mathbb{R}^m)$. We introduce $z\in \mathcal{C}^1(\mathbb{R}_+;\mathbb{R}^n)$ defined for all $t \geq 0$ by (see~\cite{artstein1982linear}):
\begin{equation}\label{eq: Artstein transformation}
z(t) = e^{D_0 A} x(t) + \int_{t-D_0}^{t} e^{(t-s)A} B u(s) \diff s .
\end{equation}
As $u = \varphi K z$, we have for all $t \geq 0$, 
\begin{align}
\dot{z}(t) 
& = (A+\varphi(t)BK) z(t) \label{eq: proof thm 1 - ODE z for t geq 0} \\
& \phantom{=}\, + e^{D_0 A}BK \{ [\varphi z](t-D(t)) - [\varphi z](t-D_0) \} . \nonumber
\end{align} 
In particular, we have for all $t \geq t_1 \triangleq t_0 + D_0 + \delta$ that
\begin{equation}\label{eq: proof thm 1 - ODE z for t geq t0+D0+delta}
\dot{z}(t) = A_\mathrm{cl} z(t) + e^{D_0 A}BK \{ z(t-D(t)) - z(t-D_0) \}
\end{equation} 
with $A_\mathrm{cl} = A + BK$ Hurwitz and the continuously differentiable initial condition $\left. z \right\vert_{[t_0,t_1]}$. Applying Lemma~\ref{lem: lemma 1}, we obtain that $\Vert z(t) \Vert \leq C_0 e^{-\kappa (t-t_1)} \Vert z(t_1 + \cdot) \Vert_W$ for $t \geq t_1$.

We introduce $V(t) = \Vert z(t) \Vert^2 / 2$ for $t \geq 0$. The use of the Young's inequality shows that there exist constants $\gamma_1,\gamma_2 > 0$, independent of $x_0$ and $D$, such that for all $t \geq 0$,
$ 
\dot{V}(t) 
\leq \gamma_1 V(t) + \gamma_2 [\varphi(t-D(t))]^2 V(t-D(t))
\phantom{\leq}\, + \gamma_2 [\varphi(t-D_0)]^2 V(t-D_0)
$.
We show by induction that, for any $n \in \mathbb{N}^*$, there exists a constant $c_n > 0$, independent of $x_0$ and $D$, such that $V(t) \leq c_n^2 \Vert x_0 \Vert^2 /2$ for all $t \in [0,n(D_0 - \delta)]$. In the case $n=1$, we have for all $t \in [0,D_0-\delta]$, $\varphi(t-D(t)) = \varphi(t-D_0) = 0$. Thus $\dot{V}(t) \leq \gamma_1 V(t)$ and we obtain that the property holds true with $c_1 = e^{\gamma_1 (D_0 - \delta)/2} \Vert e^{D_0 A} \Vert$. Assume that $V(t) \leq c_n^2 \Vert x_0 \Vert^2 /2$ for all $t \in [0,n(D_0 - \delta)]$. Then, for all $t \in [0,(n+1)(D_0 - \delta)]$, we have $t-D(t) \leq n(D_0 - \delta)$ and $t-D_0 \leq n(D_0 - \delta)$, yielding $\dot{V}(t) \leq \gamma_1 V(t) + \gamma_2 c_n^2 \Vert x_0 \Vert^2$. A straightforward integration shows the existence of the claimed $c_{n+1} > 0$.

Let $n_0 \geq 1$ be such that $n_0 (D_0 - \delta) \geq t_1$. This yields $\sup\limits_{t \in [0,t_1]} \Vert z(t) \Vert \leq c_{n_0} \Vert x_0 \Vert$. From (\ref{eq: proof thm 1 - ODE z for t geq 0}), we infer the existence of a constant $\tilde{c}_0 > 0$, independent of $x_0$ and $D$, such that $\sup\limits_{t \in [0,t_1]} \Vert \dot{z}(t) \Vert \leq \tilde{c}_0 \Vert x_0 \Vert$. From the definition of $\Vert \cdot \Vert_W$, we obtain that $\Vert z(t_1 + \cdot) \Vert_W \leq \tilde{c}_1 \Vert x_0 \Vert$ with $\tilde{c}_1 = \sqrt{c_{n_0}^2 + (D_0+\delta) \tilde{c}_0^2}$. We deduce that $\Vert z(t) \Vert \leq \tilde{C}_0 e^{-\kappa t} \Vert x_0 \Vert$ for all $t \geq 0$ with $\tilde{C}_0 =  e^{\kappa t_1} \max( C_0 \tilde{c}_1 , c_{n_0} ) > 0$. The conclusion follows from straightforward estimations of $u = \varphi K z$ and (\ref{eq: Artstein transformation}). \qed

\begin{rem}
In Theorem~\ref{thm: exp stability for finite-dimensional systems}, the initial control input is identically zero, i.e., $u_0 \triangleq \left. u \right\vert_{[-D_0-\delta,0]} = 0$. This can be obtained in practice by initially applying a zero control input. This avoids the necessity of 1) regularity assumptions on $u_0$; 2) the introduction of compatibility conditions restricting the admissible initial conditions $x_0$ (see Theorem~\ref{thm: exp stab - small gain}); 3) the explicit knowledge of $u_0$ to initialize the computation of the predictor feedback. Note that in the case of an actuator exhibiting a dynamical behavior, the initial actuator state is, in general, non zero. In this case, one could augment the state of the plant with the dynamics of the actuator. In this setting, the initial condition of the actuator is captured by $x_0$.
\end{rem}

\subsection{Applications}\label{subsec: applications}

Using the LMI solvers of \textsc{Matlab}~R2017b, we compare the application of the results of: (T1) Theorem~\ref{thm: exp stab - small gain} from~\cite{karafyllis2013delay}; (T2) the LMI condition from~\cite[Thm~2]{selivanov2016predictor}; (T3) Theorem~\ref{thm: exp stability for finite-dimensional systems}. The examples are extracted from~\cite{li2014robustness}. 

\begin{exmp}
With the matrices
\begin{equation*}
A = \begin{bmatrix} 0 & 1 \\ -1 & 1 \end{bmatrix}, \quad
B = \begin{bmatrix} 0 \\ 1 \end{bmatrix}, \quad
K = \begin{bmatrix} -1 & -3 \end{bmatrix}, \quad
\end{equation*}
the closed-loop poles are located in $-1 \pm j$. For $D_0 = 1\,\mathrm{s}$, we obtain (T1) $\delta = 0.0212$ ($\delta_\mathrm{E} = 0.0400$); (T2) with $\kappa = 0.2$, $\delta = 0.0563$; (T3) with $\kappa = 0.2$, $\delta = 0.0780$.
\end{exmp}

\begin{exmp}
With the matrices 
\begin{equation*}
A = \begin{bmatrix} -2/3 & -1 & 5/3 \\ 0 & -1 & 0 \\ 1/3 & -1 & 2/3 \end{bmatrix}, \quad
B = \begin{bmatrix} 1 & -1 \\ 0 & 2 \\ -2 & 1 \end{bmatrix},
\end{equation*}
\begin{equation*}
K = \begin{bmatrix} 0.3572 & -0.4853 & 1.1281 \\ 0.3925 & -0.5660 & 0.4235 \end{bmatrix},
\end{equation*}
the closed-loop poles are located in $-1 \pm j$ and $-2$. For $D_0 = 1\,\mathrm{s}$, we obtain (T1) $\delta = 0.0147$ ($\delta_\mathrm{E} = 0.0391$); (T2) with $\kappa = 0.2$, $\delta = 0.0591$; (T3) with $\kappa = 0.2$, $\delta = 0.0796$.
\end{exmp}

\section{Extension to the feedback stabilization of a class of diagonal infinite-dimensional systems}\label{sec: infinite dimensional systems}

We extend the results of Theorem~\ref{thm: exp stability for finite-dimensional systems} to the feedback stabilization of a class of diagonal (infinite-dimensional) boundary control systems exhibiting a finite number of unstable modes by means of a boundary control input that is subject to an uncertain and time-varying delay. In the sequel,  $(\mathcal{H},\left<\cdot,\cdot\right>_\mathcal{H})$ is a separable $\mathbb{K}$-Hilbert space.

\subsection{Problem setting}

Let $D_0>0$ and $\delta \in (0,D_0)$ be given. We consider the abstract boundary control system~\cite{Curtain2012}:
\begin{equation}\label{def: boundary control system}
\left\{\begin{split}
\dfrac{\mathrm{d} X}{\mathrm{d} t}(t) & = \mathcal{A} X(t) , & t \geq 0 \\
\mathcal{B} X (t) & = \tilde{u}(t) \triangleq u(t-D(t)) , & t \geq 0 \\
X(0) & = X_0 
\end{split}\right.
\end{equation}
with
\begin{itemize}
\item $\mathcal{A} : D(\mathcal{A}) \subset \mathcal{H} \rightarrow \mathcal{H}$ a linear (unbounded) operator;
\item $\mathcal{B} : D(\mathcal{B}) \subset \mathcal{H} \rightarrow \mathbb{K}^m$ with $D(\mathcal{A}) \subset D(\mathcal{B})$ a linear boundary operator;
\item $u : [-D_0-\delta,+\infty) \rightarrow \mathbb{K}^m$ with $\left. u \right\vert_{[-D_0-\delta,0)} = 0$ the boundary control;
\item $D:\mathbb{R}_+ \rightarrow [D_0-\delta,D_0+\delta]$ a time-varying delay.
\end{itemize}
It is assumed that $(\mathcal{A},\mathcal{B})$ is a boundary control system:
\begin{enumerate}
\item the disturbance-free operator $\mathcal{A}_0$, defined on the domain $D(\mathcal{A}_0) \triangleq D(\mathcal{A}) \cap \mathrm{ker}(\mathcal{B})$ by $\mathcal{A}_0 \triangleq \left.\mathcal{A}\right|_{D(\mathcal{A}_0)}$, is the generator of a $C_0$-semigroup $S$ on $\mathcal{H}$;
\item there exists a bounded operator $B \in \mathcal{L}(\mathbb{K}^m,\mathcal{H})$, called a lifting operator, such that $\mathrm{R}(B) \subset D(\mathcal{A})$, $\mathcal{A}B \in \mathcal{L}(\mathbb{K}^m,\mathcal{H})$, and $\mathcal{B}B = I_{\mathbb{K}^m}$;
\end{enumerate}
where $\mathrm{ker}(\mathcal{B})$ stands for the kernel of $\mathcal{B}$ and $\mathrm{R}(B)$ denotes the range of $B$. 

In the following developments, we assume that the boundary control system exhibits a diagonal structure:
\begin{assum}\label{assum: A1}
The disturbance-free operator $\mathcal{A}_0$ is a Riesz spectral operator~\cite{Curtain2012}, i.e., is a linear and closed operator with simple eigenvalues $\lambda_n$ and corresponding eigenvectors $\phi_n \in D(\mathcal{A}_0)$, $n \in \mathbb{N}^*$, that satisfy:
\begin{enumerate}
\item $\left\{ \phi_n , \; n \in \mathbb{N}^* \right\}$ is a Riesz basis~\cite{christensen2016introduction}:
\begin{enumerate}
\item $\overline{ \underset{n\in\mathbb{N}^*}{\mathrm{span}_\mathbb{K}} \;\phi_n } = \mathcal{H}$;
\item there exist constants $m_R, M_R \in \mathbb{R}_+^*$ such that for all $N \in \mathbb{N}^*$ and all $\alpha_1 , \ldots , \alpha_N \in \mathbb{K}$,
\begin{equation}\label{eq: Riesz basis - inequality}
m_R \sum\limits_{n=1}^{N} \vert \alpha_n \vert^2
\leq
\left\Vert \sum\limits_{n=1}^{N} \alpha_n \phi_n \right\Vert_\mathcal{H}^2
\leq
M_R \sum\limits_{n=1}^{N} \vert \alpha_n \vert^2 .
\end{equation}
\end{enumerate}
\item The closure of $\{ \lambda_n , \; n \in \mathbb{N}^* \}$ is totally disconnected, i.e. for any distinct $a,b \in \overline{ \{ \lambda_n , \; n \in \mathbb{N}^* \} }$, $[a,b] \not\subset \overline{ \{ \lambda_n , \; n \in \mathbb{N}^* \} }$.
\end{enumerate}
\end{assum}

We also assume that the system presents a finite number of unstable modes and that the set composed of the real part of the stable modes does not accumulate at 0:
\begin{assum}\label{assum: A2}
There exist $N_0 \in \mathbb{N}^*$ and $\alpha \in \mathbb{R}_+^*$ such that $\operatorname{Re} \lambda_n \leq - \alpha$ for all $n \geq N_0 + 1$. 
\end{assum}

As $\left\{ \phi_n , \; n \in \mathbb{N}^* \right\}$ is a Riesz basis, we can introduce its biorthogonal sequence $\left\{ \psi_n , \; n \in \mathbb{N}^* \right\}$, i.e., $\left< \phi_k , \psi_l \right>_\mathcal{H} = \delta_{k,l} \in\{0,1\}$ with $\delta_{k,l}=1$ if and only if $k=l$. Then, we have for all $x \in \mathcal{H}$ the following series expansion: $x = \sum\limits_{n \geq 1} \left< x , \psi_n \right>_\mathcal{H} \phi_n$. As $\mathcal{A}_0$ is a Riesz-spectral operator, then $\psi_n$ is an eigenvector of the adjoint operator $\mathcal{A}_0^*$ associated with the eigenvalue $\overline{\lambda_n}$. 

\subsection{Spectral decomposition and finite dimensional truncated model}
Under the assumption that $\tilde{u} \in \mathcal{C}^2([0,+\infty);\mathbb{K}^m)$ and $X_0 \in D(\mathcal{A})$ such that $\mathcal{B}X_0 = \tilde{u}(0) = u(0-D(0)) = 0$ (i.e., $X_0 \in D(\mathcal{A}_0)$), there exists a unique classical solution $X \in \mathcal{C}^0(\mathbb{R}_+;D(\mathcal{A})) \cap \mathcal{C}^1(\mathbb{R}_+;\mathcal{H})$ of (\ref{def: boundary control system}); see, e.g., \cite[Th.~3.3.3]{Curtain2012}. Then, 
\begin{equation*}
X(t) 
= \sum\limits_{n \in \mathbb{N}^*} \left< X(t) , \psi_n \right>_\mathcal{H} \phi_n 
= \sum\limits_{n \in \mathbb{N}^*} c_n(t) \phi_n ,
\end{equation*}
where $c_n(t) \triangleq \left< X(t) , \psi_n \right>_\mathcal{H}$. We infer that $c_n \in \mathcal{C}^1(\mathbb{R}_+;\mathbb{K})$ and, from (\ref{def: boundary control system}), we have for all $t \geq 0$ the following spectral decomposition~\cite{lhachemi2018iss}:
\begin{align}
& \dot{c}_n(t) \nonumber \\
& = \left< \mathcal{A} X(t) , \psi_n \right>_\mathcal{H} \nonumber \\
& = \left< \mathcal{A}_0 \left\{ X(t) - B \tilde{u}(t) \right\} , \psi_n \right>_\mathcal{H} + \left< \mathcal{A} B \tilde{u}(t) , \psi_n \right>_\mathcal{H} \nonumber \\
& = \left< X(t) - B \tilde{u}(t) , \mathcal{A}_0^* \psi_n \right>_\mathcal{H} + \left< \mathcal{A} B \tilde{u}(t) , \psi_n \right>_\mathcal{H} \nonumber \\
& = \lambda_n c_n(t) - \lambda_n \left< B \tilde{u}(t) , \psi_n \right>_\mathcal{H} + \left< \mathcal{A} B \tilde{u}(t) , \psi_n \right>_\mathcal{H} , \label{eq: coeff in Riesz basis ODE}
\end{align}
where it has been used that $\mathcal{B} \left\{ X(t) - B \tilde{u}(t) \right\} = \tilde{u}(t) - \tilde{u}(t) = 0$, showing that $X(t) - B \tilde{u}(t) \in D(\mathcal{A}) \cap \mathrm{ker}(\mathcal{B}) = D(\mathcal{A}_0)$.

Let $\mathcal{E} = (e_1,e_2,\ldots,e_m)$ be the canonical basis of $\mathbb{K}^m$. Introducing $b_{n,k} \triangleq - \lambda_n \left< B e_k , \psi_n \right>_\mathcal{H} + \left< \mathcal{A} B e_k , \psi_n \right>_\mathcal{H}$, we obtain from (\ref{eq: coeff in Riesz basis ODE}) that the following linear ODE holds true for all $t \geq 0$
\begin{equation}\label{eq: ODE satisfies by Y}
\dot{Y}(t) 
= A_{N_0} Y(t) + B_{N_0} u(t-D(t)) ,
\end{equation}
where $A_{N_0} = \mathrm{diag}(\lambda_1,\ldots,\lambda_{N_0}) \in \mathbb{K}^{N_0 \times N_0}$, $B_{N_0} = (b_{n,k})_{1 \leq n \leq N_0 , 1 \leq k \leq m} \in \mathbb{K}^{N_0 \times m}$, and
\begin{equation}\label{eq: def Y}
Y(t) = 
\begin{bmatrix}
c_1(t) \\ \vdots \\ c_{N_0}(t)
\end{bmatrix}
=
\begin{bmatrix}
\left< X(t) , \psi_1 \right>_\mathcal{H} \\ \vdots \\ \left< X(t) , \psi_{N_0} \right>_\mathcal{H}
\end{bmatrix} 
\in \mathbb{K}^{N_0} .
\end{equation}

Under the following assumption, we obtain the existence of a feedback gain $\mathbb{K}^{m \times N_0}$ such that $A_{N_0} + B_{N_0} K$ is Hurwitz.
\begin{assum}\label{assum: A3}
$(A_{N_0},B_{N_0})$ is stabilizable.
\end{assum}

Then, we can employ the strategy presented in Section~\ref{sec: finite dimensional systems} to ensure the exponential feedback stabilization of the finite-dimensional truncated dynamics (\ref{eq: ODE satisfies by Y}). The objective is now to assess that such a strategy ensures the stabilization of the full infinite-dimensional system.

\begin{rem}
In general, even for problems originally defined over the real field $\mathbb{K} = \mathbb{R}$, the spectral decomposition (\ref{eq: ODE satisfies by Y}) might be complex-valued due to the incursion into the complex plan to define the eigenstructures $(\lambda_n,\phi_n)$ of the system ; typical examples of such systems are strings and beams. Consequently, we need in the sequel the following complex-version of Theorem~\ref{thm: exp stability for finite-dimensional systems}.
\end{rem}

\begin{cor}\label{cor: exp stability for finite-dimensional systems - complex case}
In the context of Theorem~\ref{thm: exp stability for finite-dimensional systems} but with complex-valued $A$, $B$, $K$, and $x_0$, i.e., $A \in \mathbb{C}^{n \times n}$, $B \in \mathbb{C}^{n \times m}$, $K \in \mathbb{C}^{m \times n}$, and $x_0 \in \mathbb{C}^{n}$, the conclusions of Theorem~\ref{thm: exp stability for finite-dimensional systems} hold true with $M = \mathcal{R}(A_\mathrm{cl})$ and $N = \mathcal{R}(e^{D_0 A} BK)$. In this case, the matrices of (\ref{eq: preleminary lemma LMI condition}) are such that $P_1,Q\in\mathbb{S}_{2n}^{+*}$ and $P_2, P_3 \in \mathbb{R}^{2n \times 2n}$.
\end{cor}

\textbf{Proof.} For $z(t) \in \mathbb{C}^n$, we infer that (\ref{eq: proof thm 1 - ODE z for t geq t0+D0+delta}) is equivalent to
\begin{align*}
\dot{\tilde{z}}(t)
& = 
\mathcal{R}(A_\mathrm{cl}) 
\tilde{z}(t) \\
& \phantom{=}\, + \mathcal{R}(e^{D_0 A} BK)
\left\{
\tilde{z}(t-D(t))
-
\tilde{z}(t-D_0)
\right\}
\end{align*}
with $\tilde{z}(t) = \begin{bmatrix} \operatorname{Re} z(t)^\top & \operatorname{Im} z(t)^\top \end{bmatrix}^\top \in \mathbb{R}^{2n}$. Furthermore, as $A_\mathrm{cl}$ is assumed Hurwitz, so is $\mathcal{R}(A_\mathrm{cl})$. Then, the conclusion follows from the proof of Theorem~\ref{thm: exp stability for finite-dimensional systems}. \qed

\subsection{Dynamics of the closed-loop system}
Let $D_0,t_0 > 0$ and $\delta \in (0,D_0)$ be given. Let $\varphi \in \mathcal{C}^2(\mathbb{R};\mathbb{R})$ be a transition signal over $[0,t_0]$ and $D \in \mathcal{C}^2(\mathbb{R}_+;\mathbb{R})$ be a time-varying delay such that $\vert D - D_0 \vert \leq \delta$. The dynamics of the closed-loop system takes the form (see~\cite{lhachemi2018feedback} for the nominal case $D(t)=D_0$):
\begin{subequations}
\begin{align}
\dfrac{\mathrm{d} X}{\mathrm{d} t}(t) & = \mathcal{A} X(t) , \label{def: boundary control system - closed-loop - first}\\
\mathcal{B} X (t) & = \tilde{u}(t) = u(t-D(t)) , \\
\left. u \right\vert_{[-D_0-\delta,0]} & = 0 , \\
u(t) & = \varphi(t) K e^{D_0 A_{N_0}} Y(t) \label{def: boundary control system - closed-loop - control input} \\
& \phantom{=} + \varphi(t) K \int_{t-D_0}^{t} e^{(t-s)A_{N_0}} B_{N_0} u(s) \diff s , \nonumber \\
X(0) & = X_0 \label{def: boundary control system - closed-loop - last}
\end{align}
\end{subequations}
for any $t \geq 0$ with $Y$ given by (\ref{eq: def Y}). The gain $K \in \mathbb{K}^{m \times N_0}$ is selected such that $A_\mathrm{cl} \triangleq A_{N_0} + B_{N_0} K$ is Hurwitz.

\begin{lem}\label{lemma: well-posedness closed-loop PDE}
Let $(\mathcal{A},\mathcal{B})$ be an abstract boundary control system such that Assumptions~\ref{assum: A1}, \ref{assum: A2}, and~\ref{assum: A3} hold true. For any $X_0 \in D(\mathcal{A}_0)$ and  $D \in \mathcal{C}^2(\mathbb{R}_+;\mathbb{R})$ such that $\vert D - D_0 \vert \leq \delta < D_0$, the closed-loop system (\ref{def: boundary control system - closed-loop - first}-\ref{def: boundary control system - closed-loop - last}) admits a unique classical solution $X \in \mathcal{C}^0(\mathbb{R}_+;D(\mathcal{A})) \cap \mathcal{C}^1(\mathbb{R}_+;\mathcal{H})$. The associated control law $u$ is the unique solution of the implicit equation (\ref{def: boundary control system - closed-loop - control input}) and is of class $\mathcal{C}^2([-D_0-\delta,+\infty);\mathbb{K}^m)$. It can be written under the form $u = \varphi K Z$ with, for all $t \geq 0$,
\begin{equation}\label{eq: Artstein transform}
Z(t) \triangleq e^{D_0 A_{N_0}} Y(t) + \int_{t-D_0}^{t} e^{(t-s)A_{N_0}} B_{N_0} u(s) \diff s ,
\end{equation}
which is such that $Z \in \mathcal{C}^2(\mathbb{R}_+;\mathbb{K}^{N_0})$ with for all $t \geq 0$,
\begin{align}
\dot{Z}(t)
& = (A_{N_0} + \varphi(t) B_{N_0} K) Z(t) \label{eq: IDS ODE Z} \\
& \phantom{=}\, + e^{D_0 A_{N_0}} B_{N_0}K  \{ [\varphi Z](t-D(t)) - [\varphi Z](t-D_0) \} . \nonumber
\end{align}
\end{lem}

The proof of Lemma~\ref{lemma: well-posedness closed-loop PDE} relies on the invertibility of the Artstein transformation~\cite{bresch2018new} and on the fact that, for any $t \in [n(D_0-\delta),(n+1)(D_0-\delta)]$ with $n \in \mathbb{N}$, the actual control input $\tilde{u}(t)$ is such that $\tilde{u}(t)=0$ for $n=0$ and depends only on the system state $X$ (via $Y$) over the range of time $[0,n(D_0-\delta)]$ when $n \geq 1$. Therefore, the existence of a classical solution $X$ for the closed-loop system (\ref{def: boundary control system - closed-loop - first}-\ref{def: boundary control system - closed-loop - last}) can be shown by induction using classical results on boundary control systems with boundary input of class $\mathcal{C}^2$ (see, e.g., \cite[Th.~3.3.3]{Curtain2012}). Such a regularity of the control input follows first from the fact that the control law $u$ implicitly defined by (\ref{def: boundary control system - closed-loop - control input}) via the Artstein transformation is of class $\mathcal{C}^0$ (see~\cite{bresch2018new}) and then from (\ref{eq: Artstein transform}-\ref{eq: IDS ODE Z}). A detailed proof in the case $\delta = 0$, i.e., $D(t)=D_0$, can be found in~\cite[Section~IV.B]{lhachemi2018feedback} and, based on the above remarks, can be extended in a straightforward manner to the configuration of Lemma~\ref{lemma: well-posedness closed-loop PDE}.

\subsection{Exponential stability of the closed-loop system}

The exponential stability of the closed-loop system (\ref{def: boundary control system - closed-loop - first}-\ref{def: boundary control system - closed-loop - last}) in the nominal case $D(t)=D_0$ has been assessed in~\cite{lhachemi2018feedback}. The contribution of this paper relies on the following robustness assessment of the control strategy with respect to uncertain and time-varying delays $D(t)$.

\begin{thm}\label{thm: IDS - main result}
Let $(\mathcal{A},\mathcal{B})$ be an abstract boundary control system such that Assumptions~\ref{assum: A1}, \ref{assum: A2}, and~\ref{assum: A3} hold true. There exist $\delta \in (0,D_0)$ and $\eta > 0$ such that, for any given $\delta_r > 0$, we have the existence of a constant $C_2>0$ such that, for any $X_0 \in D(\mathcal{A}_0)$ and $D \in \mathcal{C}^2(\mathbb{R}_+;\mathbb{R})$ with $\vert D - D_0 \vert \leq \delta$ and $\sup\limits_{t \in \mathbb{R}_+} \left\vert \dot{D}(t) \right\vert \leq \delta_r $, the trajectory $X$ and the control input $u$ of the closed-loop dynamics (\ref{def: boundary control system - closed-loop - first}-\ref{def: boundary control system - closed-loop - last}) satisfy $\Vert X(t) \Vert_\mathcal{H} + \Vert u(t) \Vert \leq C_2 e^{-\eta t} \Vert X_0 \Vert_\mathcal{H}$ for all $t \geq 0$. In particular, this conclusion holds true for any $\delta \in (0,D_0)$ such that $\Theta(\delta,0) \prec 0$ is feasible with 
\begin{itemize}
\item in the case $\mathbb{K}=\mathbb{R}$, $M = A_{N_0} + B_{N_0} K$, $N = e^{D_0 A_{N_0}} B_{N_0}K$, $P_1,Q\in\mathbb{S}_n^{+*}$, and $P_2, P_3 \in \mathbb{R}^{n \times n}$;
\item in the case $\mathbb{K}=\mathbb{C}$, $M = \mathcal{R}(A_{N_0} + B_{N_0} K)$, $N = \mathcal{R}(e^{D_0 A_{N_0}} B_{N_0}K)$, $P_1,Q\in\mathbb{S}_{2n}^{+*}$, and $P_2, P_3 \in \mathbb{R}^{2n \times 2n}$.
\end{itemize}
Furthermore, if $\kappa > 0$ is such that $\Theta(\delta,\kappa) \preceq 0$ is feasible, then the decay rate $\eta$ can be taken as any element of $(0,\kappa]$ if $\alpha > \kappa$ or $(0,\alpha)$ if $\alpha \leq \kappa$.
\end{thm}

\textbf{Proof.}
Let $\delta \in (0,D_0)$ and $\kappa > 0$ be such that $\Theta(\delta,\kappa) \preceq 0$ is feasible (see Lemma~\ref{lem: lemma 2}). We introduce $\eta \in (0,\kappa]$ if $\alpha > \kappa$ or $\eta \in (0,\alpha)$ if $\alpha \leq \kappa$. Thus, we can select a $\epsilon \in (0,1)$ such that $\alpha_\epsilon \triangleq \alpha (1 - \epsilon) > \eta$. Let $\delta_r > 0$ be arbitrarily given. Let $X_0 \in D(\mathcal{A}_0)$ and $D \in \mathcal{C}^2(\mathbb{R}_+;\mathbb{R})$ such that $\vert D - D_0 \vert \leq \delta$ and $\sup\limits_{t \in \mathbb{R}_+} \vert \dot{D}(t) \vert \leq \delta_r$ be given. From Lemma~\ref{lemma: well-posedness closed-loop PDE}, we denote by $X \in \mathcal{C}^0(\mathbb{R}_+;D(\mathcal{A})) \cap \mathcal{C}^1(\mathbb{R}_+;\mathcal{H})$ the unique classical solution of the closed-loop system (\ref{def: boundary control system - closed-loop - first}-\ref{def: boundary control system - closed-loop - last}) and $u \in \mathcal{C}^2([-D_0-\delta,+\infty);\mathbb{K}^m)$ the associated control input. Thus (\ref{eq: ODE satisfies by Y}) holds true for all $t \geq 0$. Furthermore, as $u = \varphi K Z$ with $Z$ given by (\ref{eq: Artstein transform}) and $\left. u \right\vert_{[-D_0-\delta,0]} = 0$, we obtain from Theorem~\ref{thm: exp stability for finite-dimensional systems} that $\Vert Y(t) \Vert + \Vert u(t) \Vert \leq C_1 e^{-\kappa t} \Vert Y(0) \Vert$ and $\Vert Z(t) \Vert \leq \tilde{C}_0 e^{-\kappa t} \Vert Y(0) \Vert$ for all $t \geq 0$ with constants $C_1,\tilde{C}_0>0$ independent of $X_0$ and $D$. From (\ref{eq: Riesz basis - inequality}) and (\ref{eq: def Y}), we have that $\Vert Y(0) \Vert \leq \Vert X_0 \Vert_\mathcal{H} / \sqrt{m_R}$. This yields, along with $0 < \eta \leq \kappa$, $\Vert Y(t) \Vert + \Vert u(t) \Vert \leq C_1 e^{-\eta t} \Vert X_0 \Vert_\mathcal{H}/ \sqrt{m_R}$ and $\Vert Z(t) \Vert \leq \tilde{C}_0 e^{-\eta t} \Vert X_0 \Vert_\mathcal{H}/ \sqrt{m_R}$ for all $t \geq 0$.

In order to assess the exponential stability of the full infinite-dimensional system, we introduce for all $t \geq 0$,
\begin{equation*}
V(t) = \dfrac{1}{2} \sum\limits_{k \geq N_0 + 1} \left\vert \left< X(t) - B \tilde{u}(t) , \psi_k \right> \right\vert^2 \geq 0 ,
\end{equation*}
which is such that $V(t) \leq \Vert X(t) - B \tilde{u}(t) \Vert_\mathcal{H}^2/(2 m_R) < + \infty$ and $V \in \mathcal{C}^1(\mathbb{R}_+;\mathbb{R})$. The quantity $V(t)$ is used to derive an upper bound of $\Vert X(t) \Vert_\mathcal{H}$ as follows. Noting that
\begin{align*}
& \dfrac{1}{2} \sum\limits_{k = 1}^{N_0} \left\vert \left< X(t) - B \tilde{u}(t) , \psi_k \right> \right\vert^2 \\
& \qquad \leq \sum\limits_{k = 1}^{N_0} \left\vert \left< X(t) , \psi_k \right> \right\vert^2 + \sum\limits_{k = 1}^{N_0} \left\vert \left< B \tilde{u}(t) , \psi_k \right> \right\vert^2 \\
& \qquad \overset{(\ref{eq: def Y})}{\leq} \Vert Y(t) \Vert^2 + \sum\limits_{k \geq 1} \left\vert \left< B \tilde{u}(t) , \psi_k \right> \right\vert^2 \\
& \qquad \overset{(\ref{eq: Riesz basis - inequality})}{\leq} \Vert Y(t) \Vert^2 + \dfrac{1}{m_R} \Vert B \tilde{u}(t) \Vert_\mathcal{H}^2 ,
\end{align*}
we obtain that
\begin{align*}
& V(t) \\ 
& = \dfrac{1}{2} \sum\limits_{k \geq 1} \left\vert \left< X(t) - B \tilde{u}(t) , \psi_k \right> \right\vert^2 - \dfrac{1}{2} \sum\limits_{k = 1}^{N_0} \left\vert \left< X(t) - B \tilde{u}(t) , \psi_k \right> \right\vert^2 \\
& \overset{(\ref{eq: Riesz basis - inequality})}{\geq} \dfrac{1}{2 M_R} \Vert X(t) - B \tilde{u}(t) \Vert_\mathcal{H}^2 - \Vert Y(t) \Vert^2 - \dfrac{1}{m_R} \Vert B \tilde{u}(t) \Vert_\mathcal{H}^2 .
\end{align*}
Using the triangular inequality, this yields for all $t \geq 0$,
\begin{align*}
\Vert X(t) \Vert_\mathcal{H}
& \leq
\Vert B \tilde{u}(t) \Vert_\mathcal{H} \\
& \phantom{\leq}\, + \sqrt{2M_R \left( V(t) + \Vert Y(t) \Vert^2 + \dfrac{1}{m_R} \Vert B \tilde{u}(t) \Vert_\mathcal{H}^2 \right)} . \nonumber
\end{align*}
Noting that $t-D(t) \geq t-D_0-\delta$, we have $\Vert \tilde{u}(t) \Vert = \Vert u(t-D(t)) \Vert \leq C_1 e^{\eta (D_0 + \delta)} e^{-\eta t} \Vert X_0 \Vert_\mathcal{H}/ \sqrt{m_R}$. As $B$ is bounded and $\Vert Y(t) \Vert \leq C_1 e^{-\eta t} \Vert X_0 \Vert_\mathcal{H}/ \sqrt{m_R}$, the proof will be complete if we can show the existence of $\tilde{C}_1 > 0$, independent of $X_0$ and $D$, such that $V(t) \leq \tilde{C}_1 e^{-2\eta t} \Vert X_0 \Vert_\mathcal{H}^2$. To do so, we compute for $t \geq 0$ the time derivative of $V$ as follows:
\begin{align*}
\dot{V}(t) & = \sum\limits_{k \geq N_0+1} \operatorname{Re} \left\{ 
\left< \dfrac{\mathrm{d}X}{\mathrm{d}t}(t) - B\dot{\tilde{u}}(t) , \psi_k \right>_\mathcal{H} \right. \\
& \phantom{= \sum\limits_{k \geq N_0+1} \operatorname{Re} \{}\qquad \times \left.
\overline{\left< X(t) - B\tilde{u}(t) , \psi_k \right>_\mathcal{H}} 
\right\} ,
\end{align*}
where $\dot{\tilde{u}}(t) = (1-\dot{D}(t))\dot{u}(t-D(t))$. Using (\ref{eq: coeff in Riesz basis ODE}), Assumption~\ref{assum: A2}, and the Young Inequality (Y.I.), we obtain that
\begin{align}
\dot{V}(t) 
& \overset{(\ref{eq: coeff in Riesz basis ODE})}{=} \sum\limits_{k \geq N_0+1} \operatorname{Re}(\lambda_k) \left\vert \left< X(t) - B\tilde{u}(t) , \psi_k \right>_\mathcal{H} \right\vert^2 \nonumber \\
& \phantom{\overset{(\ref{eq: coeff in Riesz basis ODE})}{=}}\, + \sum\limits_{k \geq N_0+1} \operatorname{Re} \left\{ 
\left( \left< \mathcal{A}B\tilde{u}(t) , \psi_k \right>_\mathcal{H} - \left< B\dot{\tilde{u}}(t) , \psi_k \right>_\mathcal{H} \right) \right. \nonumber \\
& \phantom{= \sum\limits_{k \geq N_0+1} \operatorname{Re} \{}\qquad\qquad\qquad \times \left.
\overline{\left< X(t) - B\tilde{u}(t) , \psi_k \right>_\mathcal{H}} 
\right\} \nonumber \\
& \leq - 2 \alpha V(t) \nonumber \\
& \phantom{\leq}\, + \sum\limits_{k \geq N_0+1} 
\left( \left\vert \left< \mathcal{A}B\tilde{u}(t) , \psi_k \right>_\mathcal{H} \right\vert + \left\vert \left< B\dot{\tilde{u}}(t) , \psi_k \right>_\mathcal{H} \right\vert \right) \nonumber \\
& \phantom{= \sum\limits_{k \geq N_0+1}}\qquad\qquad\qquad \times
\left\vert \left< X(t) - B\tilde{u}(t) , \psi_k \right>_\mathcal{H} \right\vert \nonumber \\
& \overset{\mathrm{Y.I.}}{\leq} - 2 \alpha_\epsilon V(t) \label{eq: preliminary estimation dot_V} \\
& \phantom{\overset{\mathrm{Y.I.}}{\leq}}\, + \dfrac{1}{2\epsilon\alpha} \sum\limits_{k \geq N_0+1} 
\left( \left\vert \left< \mathcal{A}B\tilde{u}(t) , \psi_k \right>_\mathcal{H} \right\vert^2 + \left\vert \left< B\dot{\tilde{u}}(t) , \psi_k \right>_\mathcal{H} \right\vert^2 \right) . \nonumber 
\end{align}
For all $t \geq D_0+\delta+t_0$, as $t-D(t) \geq t_0$ and thus $\varphi(t-D(t))=1$, we have
\begin{equation*}
\tilde{u}(t) = u(t-D(t)) = K Z(t-D(t))
= \sum\limits_{i=1}^{m} \left\{ K_i Z(t-D(t)) \right\} e_i
\end{equation*}
where $K_i$ stands for $i$-th line of $K$. We deduce that
\begin{align*}
& \sum\limits_{k \geq N_0+1} 
\left\vert \left< \mathcal{A}B\tilde{u}(t) , \psi_k \right>_\mathcal{H} \right\vert^2 \\
& \qquad \leq m \sum\limits_{i=1}^m \sum\limits_{k \geq 1} \vert \left< \mathcal{A}Be_i , \psi_k \right> \vert^2 \vert K_i Z(t-D(t)) \vert^2 \\
& \qquad \overset{(\ref{eq: Riesz basis - inequality})}{\leq} \dfrac{m}{m_R} \left( \sum\limits_{i=1}^m \Vert \mathcal{A}Be_i \Vert_\mathcal{H}^2 \Vert K_i \Vert^2 \right) \Vert Z(t-D(t)) \Vert^2 .
\end{align*}
Similarly, for all $t \geq t_1 \triangleq 2(D_0+\delta)+t_0$,
\begin{align*}
& \dot{u}(t-D(t)) \\
& = K \dot{Z}(t-D(t)) \\
& \overset{(\ref{eq: IDS ODE Z})}{=} K A_\mathrm{cl} Z(t-D(t)) \\
& \phantom{\overset{(\ref{eq: IDS ODE Z})}{=}}\, + K \tilde{B}_{N_0} K \{ Z(t-2D(t)) - Z(t-D(t)-D_0) \}  \\
& = \sum\limits_{i=1}^{m} \left\{ K_i A_\mathrm{cl} Z(t-D(t)) \right\} e_i \\
& \phantom{=}\, + \sum\limits_{i=1}^{m} \left\{ K_i \tilde{B}_{N_0} K \{ Z(t-2D(t)) - Z(t-D(t)-D_0) \} \right\} e_i 
\end{align*}
with $\tilde{B}_{N_0} \triangleq e^{D_0 A_{N_0}} B_{N_0}$. We deduce that
\begin{align*}
& \sum\limits_{k \geq N_0+1} 
\left\vert \left< B\dot{\tilde{u}}(t) , \psi_k \right>_\mathcal{H} \right\vert^2 \\
& \leq \dfrac{2\beta m}{m_R} \left( \sum\limits_{i=1}^m \Vert Be_i \Vert_\mathcal{H}^2 \Vert K_i A_\mathrm{cl} \Vert^2 \right) \Vert Z(t-D(t)) \Vert^2 \\
& \phantom{\leq}\, + \dfrac{2\beta m}{m_R} \left( \sum\limits_{i=1}^m \Vert Be_i \Vert_\mathcal{H}^2 \Vert K_i \tilde{B}_{N_0} K \Vert^2 \right) \\
& \phantom{\leq\, +}\, \qquad\qquad \times \Vert Z(t-2D(t)) - Z(t-D(t)-D_0) \Vert^2 ,
\end{align*}
where $\beta \triangleq \left( 1 + \delta_r \right)^2$.
Thus, introducing the constants $k_1,k_2 \geq 0$ defined by:
\begin{align*}
k_1 & = \dfrac{m}{2\epsilon\alpha m_R} \sum\limits_{i=1}^m \left\{ \Vert \mathcal{A}Be_i \Vert_\mathcal{H}^2 \Vert K_i \Vert^2 + 2\beta \Vert Be_i \Vert_\mathcal{H}^2 \Vert K_i A_\mathrm{cl} \Vert^2 \right\} , \\
k_2 & = \dfrac{\beta m}{\epsilon\alpha m_R}  \sum\limits_{i=1}^m \Vert B e_i \Vert_\mathcal{H}^2 \Vert K_i \tilde{B}_{N_0} K \Vert^2,
\end{align*}
we obtain that, for all $t \geq t_1$, $\dot{V}(t) \leq - 2 \alpha_\epsilon V(t) + \omega(t)$ with $\omega(t) \geq 0$ defined by:
\begin{align*}
\omega(t) & = k_1 \Vert Z(t-D(t)) \Vert^2 \nonumber \\
& \phantom{=}\, + k_2 \Vert Z(t-2D(t)) - Z(t-D(t)-D_0) \Vert^2 \nonumber \\
& \leq k_3 e^{- 2 \eta t} \Vert X_0 \Vert_\mathcal{H}^2 ,
\end{align*}
where $k_3 = \tilde{C}_0^2 e^{2 \eta (D_0+\delta)} \{ k_1 + 2 k_2 e^{2 \eta D_0} (e^{2 \eta \delta} + 1) \} / m_R$. Now, as $V$ is of class $\mathcal{C}^1$ over $\mathbb{R}_+$, we obtain after integration that, for all $t \geq t_1$,
\begin{align*}
V(t) 
& \leq e^{-2 \alpha_\epsilon(t-t_1)} V(t_1) 
+ e^{-2 \alpha_\epsilon t} \int_{t_1}^t e^{2 \alpha_\epsilon \tau} \omega(\tau) \diff \tau \\
& \leq e^{-2 \eta(t-t_1)} V(t_1)
+ \dfrac{k_3}{2 (\alpha_\epsilon - \eta)} e^{-2 \eta t} \Vert X_0 \Vert_\mathcal{H}^2
\end{align*}
where it as been used that $\alpha_\epsilon > \eta$.

It remains to evaluate $V(t)$ for $0 \leq t \leq t_1$. From the exponential estimate of $Z(t)$, we have that $\Vert Z(t) \Vert \leq \tilde{C}_0 \Vert X_0 \Vert_\mathcal{H} / \sqrt{m}_R$ for all $t \geq 0$. From (\ref{eq: IDS ODE Z}), we deduce the existence of a constant $\tilde{c}_2 > 0$, independent of $X_0$ and $D$, such that $\Vert \dot{Z}(t) \Vert \leq \tilde{c}_2 \Vert X_0 \Vert_\mathcal{H}$ for all $t \geq 0$. Then, from (\ref{eq: preliminary estimation dot_V}), the facts that $\sup\limits_{t \in \mathbb{R}_+} \vert \dot{D}(t) \vert \leq \delta_r$ and $\sup\limits_{t \in \mathbb{R}} \left\vert \dot{\varphi}(t) \right\vert < +\infty$, and $V(0) \leq \Vert X_0 \Vert_\mathcal{H}^2 / (2 m_R)$, we obtain the existence of a constant $\tilde{c}_3 > 0$, independent of $X_0$ and $D$, such that $V(t) \leq \tilde{c}_3 \Vert X_0 \Vert_\mathcal{H}^2$ for all $t \in [ 0 , t_1 ]$. Consequently, we obtain that $V(t) \leq \tilde{C}_1 e^{-2 \eta t} \Vert X_0 \Vert_\mathcal{H}^2$ for all $t \geq 0$ with $\tilde{C}_1 = \tilde{c}_3 e^{2 \eta t_1} + \dfrac{k_3}{2 (\alpha_\epsilon - \eta)}$. This completes the proof.\qed

\section{Illustrative example}\label{sec: illustrative example}
We consider the following one-dimensional reaction-diffusion equation on $(0,L)$ with a delayed Dirichlet boundary control:
\begin{equation*}
\left\{\begin{split}
y_{t}(t,x) & = a y_{xx}(t,x) + c y(t,x) , & (t,x)\in \mathbb{R}_+ \times (0,L) \\
\begin{bmatrix} y(t,0) \\ y(t,L) \end{bmatrix} & = u(t-D(t)) , & t > 0
\end{split}\right.
\end{equation*}
with $a,c > 0$, $y(t,x)\in\mathbb{R}$, and $u(t) \in \mathbb{R}^2$. Introducing the $\mathbb{R}$-Hilbert space $\mathcal{H} = L^2(0,L)$ with $\left< f , g \right>_\mathcal{H} = \int_0^L fg \diff x$, it is well-known that the above reaction-diffusion equation can be written under the form of the abstract boundary control system (\ref{def: boundary control system}) with $X(t) = y(t,\cdot) \in\mathcal{H}$, $\mathcal{A}f = af''+cf$ on the domain $D(\mathcal{A}) = H^2(0,L)$, and the boundary operator $\mathcal{B}f = (f(0),f(L))$ on the domain $D(\mathcal{B}) = H^1(0,L)$. An example of lifting operator $B$ associated with $(\mathcal{A},\mathcal{B})$ is given for any $(u_1,u_2) \in \mathbb{R}^2$ by $\{B(u_1,u_2)\}(x) = u_1 + (u_2 - u_1)x/L$ with $x \in (0,L)$. It is well-known that the disturbance-free operator $\mathcal{A}_0$ is a Riesz-spectral operator that generates a $C_0$-semigroup with $\lambda_n = c - a n^2 \pi^2 / L^2$ and $\phi_n (x) = \psi_n(x) = \sqrt{2/L} \sin( n \pi x / L )$, $n \geq 1$. Then, the boundary control system $(\mathcal{A},\mathcal{B})$ satisfies Assumptions~\ref{assum: A1} and~\ref{assum: A2}. Furthermore, straightforward computations show that $b_{n,1} = a n \pi \sqrt{2 / L^3}$ and $b_{n,2} = (-1)^{n+1} a n \pi \sqrt{2 / L^3}$. As the eigenvalues are simple and $b_{n,k} \neq 0$ for all $n \geq 1$ and $k \in\{1,2\}$, we obtain from the Kalman condition that $(A_{N_0},B_{N_0})$ is controllable, fulfilling Assumption~\ref{assum: A3}. Thus, one can compute a feedback gain $K \in \mathbb{K}^{m \times N_0}$ such that $A_{N_0} + B_{N_0} K$ is Hurwitz and then apply the result of Theorem~\ref{thm: IDS - main result} for ensuring the exponential stability of the closed-loop system.

\begin{figure}
\centering
\includegraphics[width=3.3in,height=1in]{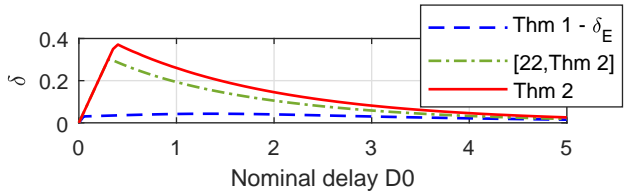}
\caption{Comparison of the estimate $\delta_\mathrm{E}$ (\ref{eq: thm exp stab - small gain condition - estimate}) derived from Theorem~\ref{thm: exp stab - small gain} with the admissible values of $\delta \in (0,D_0)$ computed for $\kappa = 0.2$ based on~\cite[Thm~2]{selivanov2016predictor} and Theorem~\ref{thm: exp stability for finite-dimensional systems}}
\label{fig: comparison small gain vs LMI}
\end{figure}

For numerical computations, we set $a=c=0.5$ and $L=2\pi$. In this configuration, we have two unstable modes $\lambda_1 = 0.375$ and $\lambda_2 = 0$ while the two first stable modes are such that $\lambda_3 = -0.625$ and $\lambda_4 = -1.5$. Setting $N_0 = 3$, the feedback gain $K \in \mathbb{R}^{2 \times 3}$ is computed to place the poles of the closed-loop truncated model $A_\mathrm{cl} = A_{N_0} + B_{N_0} K$ at $-0.75$, $-1$, and $-1.25$. Over the range $D_0 \in (0,5]$, Figure~\ref{fig: comparison small gain vs LMI} depicts: 1) the estimate $\delta_\mathrm{E}$ (\ref{eq: thm exp stab - small gain condition - estimate}) on the admissible values of $\delta > 0$ given by Theorem~\ref{thm: exp stab - small gain} taken from~\cite{karafyllis2013delay} ; 2) with decay rate $\kappa = 0.2$, the admissible values of $\delta > 0$ based on \cite[Thm~2]{selivanov2016predictor} and Theorem~\ref{thm: exp stability for finite-dimensional systems}. For the studied example, the values of $\delta$ provided by Theorem~\ref{thm: exp stability for finite-dimensional systems} are significantly less conservative. 

\begin{figure}
\centering
\includegraphics[width=3.3in,height=2.25in]{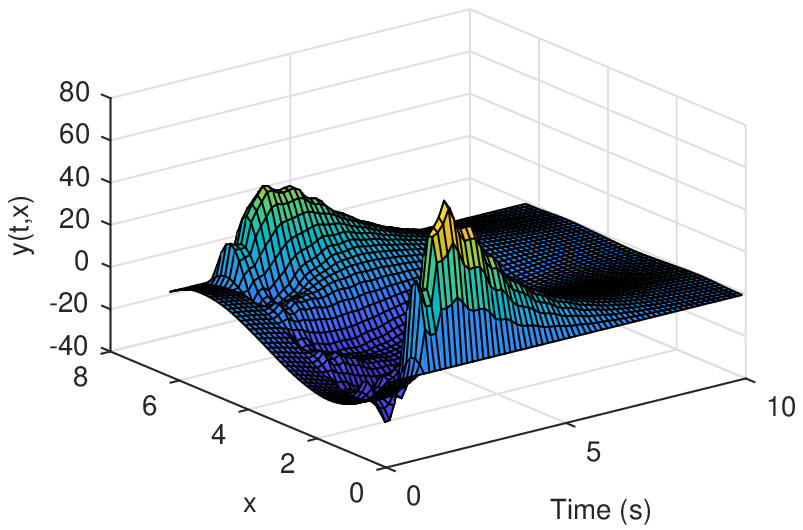}
\caption{Time evolution of $y(t)$ for the closed-loop system}
\label{fig: time evolution X}


\includegraphics[width=3.3in,height=1in]{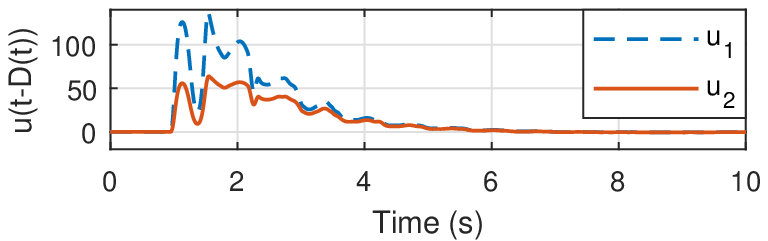}
\caption{Delayed command effort $\tilde{u}(t)=u(t-D(t))$}
\label{fig: time evolution ud}
\end{figure}

For numerical simulations, we set the nominal value of the delay to $D_0 = 1\,\mathrm{s}$. In this case, Theorem~\ref{thm: IDS - main result} ensures the exponential stability of the closed-loop system with decay rate $\kappa = 0.2$ for values of $\delta$ up to $\delta = 0.260$. We set the initial condition $X_0(x) = - x ( 2L/3 - x ) ( L - x )$ and the time-varying delay $D(t) = 1 + 0.25 \sin(3\pi t+ \pi/4)$ which is of class $\mathcal{C}^2$ and is such that $\vert D(t) - D_0 \vert \leq 0.25 \leq 0.260$ and $\vert \dot{D}(t) \vert \leq 0.75 \pi < +\infty$ for all $t \geq 0$. The transition time $t_0$ is taken as $t_0 = 0.5\,\mathrm{s}$ while the transition signal $\left.\varphi\right\vert_{[0,t_0]}$ is selected as the restriction over $[0,t_0]$ of the unique quintic polynomial function $f$ satisfying $f(0)=f'(0)=f''(0)=f'(t_0)=f''(t_0)=0$ and $f(t_0)=1$. The employed numerical scheme relies on the discretization of the reaction-diffusion equation using its first 10 modes. The time domain evolution of the closed-loop system is depicted in Figs.~\ref{fig: time evolution X}-\ref{fig: time evolution ud}. As expected from Theorem~\ref{thm: IDS - main result}, both the system state and the control input converge to zero.

\section{Conclusion}\label{sec: conclusion}
This paper discussed first the use of predictor feedback for the stabilization of finite-dimensional LTI systems in the presence of an uncertain time-varying delay in the control input. By means of a Lyapunov-Krasovskii functional, it has been derived an LMI-based sufficient condition ensuring the exponential stability of the closed-loop system for small enough variations of the time-varying delay around its nominal value. Then, this result has been extended to the feedback stabilization of a class of diagonal infinite-dimensional boundary control systems.


\bibliographystyle{plain}        
\bibliography{autosam}           





\end{document}